\documentclass{amsart}

\usepackage{amsfonts}
\usepackage{amssymb}
\usepackage{amsmath}
\usepackage{amsthm}
\usepackage{graphicx}
\usepackage{hyperref}

\theoremstyle{plain}
\newtheorem{thm}{Theorem}

\newtheorem{prop}[thm]{Proposition}

\theoremstyle{definition}

\newtheorem{remark}[thm]{Remark}

\theoremstyle{remark}

\newcommand{\R}{\mathbb{R}}

\newcommand{\C}{\mathbb{C}}
\newcommand{\La}{\Lambda}
\newcommand{\la}{\lambda}
\newcommand{\p}{\varphi}
\newcommand{\dd}{\partial}

\newcommand{\epf}{\hfill \square}
\newcommand{\sse}{\subseteq}

\newcommand{\ges}{\geqslant}
\newcommand{\x}{\times}
\newcommand{\Lag}{Lagrangian}
\newcommand{\Lags}{Lagrangians}
\newcommand{\Leg}{Legendrian}
\newcommand{\Legs}{Legendrians}
\newcommand{\Ham}{Hamiltonian}

\newcommand{\mi}{{-\infty}}
\newcommand{\std}{\text{std}}
\newcommand{\PS}{{\mathcal{PS}}}

\begin{document}

\title{Closed exact Lagrangians in the symplectization of contact manifolds}

\author[Emmy Murphy]{Emmy Murphy}
\address[Emmy Murphy]{Massachusetts Institute of Technology
  Department of Mathematics \\
  77 Massachusetts Avenue \\
  Cambridge, MA 02139, USA}
\email{e\_murphy@mit.edu}

\date{April 15, 2013}

\begin{abstract}
For a certain class of exotic contact manifolds of dimension greater than $3$, we show that there is an abundance of closed exact \Lags\ in their symplectization. All of these \Lags\ are displaceable by \Ham\ isotopy, and many of the examples are nulhomotopic.
\end{abstract}

\maketitle

Inside an exact symplectic manifold $(X, d\la)$ a \emph{formal \Lag\ embedding} is a pair $(f, \Psi_s)$, where $f:L \to X$ is a smooth embedding of a closed manifold $L$ of dimension $n = \frac12\dim X$, and $\Psi_s:TL \to TX$ is a homotopy through injective bundle homomorphisms covering $f$, so that $F_0 = df$ and $F_1$ is a map with \Lag\ image. The following theorem is a simple application of results from \cite{caps} and \cite{LPS}.

\begin{thm} \label{main thm}
Let $(Y, \ker \alpha)$ be a contact manifold of dimension $2n-1 \ges 5$. Assume that $Y$ contains a small plastikstufe with spherical core and trivial rotation. Let $(f, \Psi_s)$ be a formal \Lag\ embedding of the closed manifold $L$ into $(\R \x Y, d(e^t\alpha))$. Then $f$ is isotopic to an exact \Lag\ embedding.
\end{thm}

First we cite a theorem from \cite{LPS}.

\begin{prop}[\cite{LPS}, Theorem 1.1]\label{PS prop}
Let $(Y, \xi)$ be a contact manifold which contains a small plastikstufe with spherical core and trivial rotation, denoted by $\PS$. Let $\La \sse Y$ be a connected \Leg\ submanifold which is disjoint from $\PS$. Then $\La$ is loose.
\end{prop}

Because the specifics are not relevant to our short paper, we will omit the definition of ``small plastikstufe with spherical core and trivial rotation'', as well as the definition of loose \Legs. The plastikstufe was first defined and studied by Niederkr\"uger in \cite{orig PS}, as an obstruction to symplectic fillability. Loose \Legs\ were defined and classified up to \Leg\ isotopy in \cite{loose}. Besides the original sources, an exposition of both topics is given in \cite{LPS}. Loose \Legs\ are useful for our construction due to the following theorem from \cite{caps}; to simplify the statement we restrict our citation to the case where the ambient manifold is a symplectization.

\begin{prop}[\cite{caps}, Theorem 2.2] \label{caps prop}
Let $(Y, \ker \alpha)$ be a contact manifold of dimension $2n-1$ and let $(f, \Psi_s)$ be a formal \Lag\ embedding of the compact manifold $L$ into $([0, \infty) \x Y, d(e^t\alpha))$, so that near a collar neighborhood of the boundary $[0,\epsilon) \x \dd L$ $f$ is modeled on $f(t,p) = (t, g(p))$ for some \Leg\ embedding $g:\dd L \to Y$ and $\Psi_s$ is constant with respect to $s$. If $g(\dd L)$ is loose, then there is an isotopy of $f$ which is fixed near $\dd L$ to an exact \Lag\ embedding.
\end{prop}

Given a smooth isotopy $\p_t:X \to X$, we note that the composition $\p_1 \circ f$ can be made into a formal \Lag. This is done by taking the homotopy of bundle maps $d\p_{1-s} \circ F_s$ and translating them via a symplectic connection so that these maps cover $\p_1 \circ f$ for all $s$. This formal \Lag\ embedding is therefore non-canonical, but canonical up to homotopy.

We describe a simple model which will be useful for our construction. Inside $B^{2n} \sse \C^n$, consider the \Lag\ $L_0=\{y_i = 0, i = 1,\ldots n \}$. Notice that $\dd L_0 \sse S^{2n-1}$ is \Leg\ with respect to the standard contact form $(\frac12\sum_iy_idx_i-x_idy_i)|_{S^{2n-1}}$. We perurb $L_0$ by a \Ham\ supported near the origin so that the origin is not contained in $L_0$. Choose a radius $\gamma$ of $B^{2n}$ which is disjoint from $L_0$. Then $B^{2n} \setminus \gamma$ is symplectomorphic to $((\mi, 0] \x \R^{2n-1}, d(e^t\alpha_\std))$. Indeed, the negative Liouville flow on $B^{2n} \setminus \gamma$ is complete, and $\R^{2n-1}_\std$ is contactomorphic to $S^{2n-1}_\std$ with a single point removed. We say $L_0$ is the \emph{standard \Lag\ filling} of the \emph{standard \Leg\ unknot}, $\dd L_0$.

\noindent \emph{Proof of Theorem \ref{main thm}:} Choose a Darboux chart $U \sse Y$ which is disjoint from $\PS$. Let $\La \sse U$ be the standard \Leg\ unknot, and denote its standard \Lag\ filling by $L_0 \sse (\mi, 0] \x U$. By smooth isotopy we can assume that $f(D) = L_0$ for some disk $D \sse L$ and we can then find a homotopy of $\Psi_s$ so that on $D$, $\Psi_s = df$ for all $s \in [0,1]$. We finally arrange by smooth isotopy that $f(L) \cap (\mi, 0] \x Y = f(D)$, by finding a smooth isotopy which is fixed on $f(D)$ and moves the remainder of $f$ in the positive $t$ direction. 

Because $\La \sse U$ and $U$ is disjoint from $\PS$, Proposition \ref{PS prop} tells us that $\La$ is loose. Proposition \ref{caps prop} tells us then that $f|_{L \setminus D}$ is isotopic to an exact \Lag\ embedding. Since $f$ was already \Lag\ on $D$ this completes the construction. $\epf$

\begin{remark}
Note that any exact \Lag\ in a symplectization is \Ham\ displaceable, since any path of exact \Lag\ embeddings can be realized by \Ham\ isotopy, and in a symplectization any $t$-shift of an exact \Lag\ is again an exact \Lag.
\end{remark}

To make our result slightly more concrete, we point out that manifolds $Y$ satisfying our hypotheses are plentiful. The construction is essentially due to Etnyre and Pancholi \cite{twist}; the proposition we cite is a slight modification.

\begin{prop}[\cite{LPS}, Proposition 5.2]
Let $(Y, \xi)$ be any contact manifold of dimension larger than $3$, and let $U \sse Y$ be a Darboux chart. Then there is another contact structure $\xi^\prime$ on $Y$, so that $\xi^\prime = \xi$ outside of $U$, and $\xi^\prime$ contains a small plastikstufe with spherical core and trivial rotation inside $U$. Furthermore $\xi$ and $\xi^\prime$ are homotopic through almost contact structures.
\end{prop}

\end{document}